\pdfoutput=1
\RequirePackage{ifpdf}
\ifpdf 
\documentclass[pdftex]{sigma}
\else
\documentclass{sigma}
\fi

\numberwithin{equation}{section}

\newtheorem{Theorem}{Theorem}[section]
\newtheorem*{Theorem*}{Theorem}
\newtheorem{Corollary}[Theorem]{Corollary}

\theoremstyle{definition}

\newtheorem{Remark}[Theorem]{Remark}

\begin{document}

\allowdisplaybreaks

\newcommand{\arXivNumber}{2510.17209}

\renewcommand{\PaperNumber}{038}

\FirstPageHeading

\ShortArticleName{On Bilateral Multiple Sums and Rogers--Ramanujan Type Identities}

\ArticleName{On Bilateral Multiple Sums\\ and Rogers--Ramanujan Type Identities}

\Author{Dandan CHEN~$^{\rm ab}$ and Tianjian XU~$^{\rm a}$}

\AuthorNameForHeading{D.~Chen and T.~Xu}

\Address{$^{\rm a)}$~Department of Mathematics, Shanghai University, P.R.~China}
\EmailD{\mail{mathcdd@shu.edu.cn}, \mail{xtjmath@shu.edu.cn}}
\URLaddressD{\url{https://math.shu.edu.cn/info/1021/3135.htm}}

\Address{$^{\rm b)}$~Newtouch Center for Mathematics, Shanghai University, P.R.~China}

\ArticleDates{Received October 21, 2025, in final form March 29, 2026; Published online April 18, 2026}

\Abstract{We establish some new bilateral double-sum Rogers--Ramanujan identities involving parameters. As applications, these identities yield several new multi-sum Rogers--Ramanujan type identities. Our proofs utilize the theory of basic hypergeometric series in conjunction with the integral method.}


\Classification{11P84; 33D60}

\section{Introduction}
We begin by recalling the basic notation and definitions used throughout this paper. Assume that $|q|<1$ for convergence. The $q$-shifted factorials are defined as follows \cite{red}:
\begin{align*}
 (a;q)_0:=1,\qquad (a; q)_n := \prod_{k=0}^{n-1} \bigl(1 - a q^k\bigr),\qquad (a; q)_\infty := \prod_{k=0}^{\infty} \bigl(1 - a q^k\bigr).
\end{align*}
For negative subscripts, we define
\begin{align*}
 (a;q)_{-n}:=\frac{1}{(aq^{-n};q)_n}.
\end{align*}
We also adopt the compact product notation
\begin{align*}
 (a_1, \dots, a_m; q)_k := (a_1; q)_k \cdots (a_m; q)_k, \qquad k \in \mathbb{Z} \cup \{\infty\}.
\end{align*}\par
In 1894, Rogers \cite{Rogers-1894} discovered and proved the following fundamental identities:
\begin{align}
&\sum_{n=0}^{\infty}\frac{q^{n^2}}{(q;q)_n}=\frac{1}{\bigl(q,q^4;q^5\bigr)_{\infty}},\label{1-4-5}\\
&\sum_{n=0}^{\infty}\frac{q^{n^2+n}}{(q;q)_n}=\frac{1}{\bigl(q^2,q^3;q^5\bigr)_{\infty}}.\label{2-3-5}
\end{align}
These are now known as the Rogers--Ramanujan identities, since they were later rediscovered by Ramanujan before 1913 \cite{Hardy-1937}. These results have inspired extensive research into similar $q$-series identities, among which Slater's famous list of $130$ such identities stands out \cite{Slater-1952}. \par
A natural direction of generalization is to consider multi-sum Rogers--Ramanujan type identities. As a multi-sum generalization of \eqref{1-4-5} and \eqref{2-3-5}, the Andrews--Gordon identity \cite{Andrews-1974} stated that for integers $k>1$ and $1\leq i\leq k$,
\begin{align}\label{a-g-identity}
\sum_{n_1, \dots, n_{k-1} \geq 0} \frac{q^{N_1^2 + \dots + N_{k-1}^2 + N_i + \dots + N_{k-1}}}{(q; q)_{n_1} \cdots (q; q)_{n_{k-2}} (q; q)_{n_{k-1}}} = \frac{\bigl(q^i, q^{2k + 1 - i}, q^{2k + 1}; q^{2k + 1}\bigr)_{\infty}}{(q; q)_{\infty}},
\end{align}
where $N_j=n_j+\cdots+n_{k-1}$ for $1\leq j\leq k-1$ and $N_k=0$. Bressoud \cite{Bressoud-1980} later provided an even-modulus analogue: for integers $k>1$ and $1\leq i\leq k$,
\begin{align}\label{b-identity}
\sum_{n_1, \dots, n_{k-1} \geq 0} \frac{q^{N_1^2 + \dots + N_{k-1}^2 + N_i + \dots + N_{k-1}}}{(q; q)_{n_1} \cdots (q; q)_{n_{k-2}} \bigl(q^2; q^2\bigr)_{n_{k-1}}} = \frac{\bigl(q^i, q^{2k - i}, q^{2k }; q^{2k}\bigr)_{\infty}}{(q; q)_{\infty}},
\end{align}
with $N_j$ defined as above.\par
Rogers--Ramanujan type identities play a significant role in several fields, including combinatorics, mathematical physics and number theory. In particular, they reveal profound connections between the theory of $q$-series and modular forms. For instance, after multiplying with suitable powers of $q$, the right-hand sides of identities \eqref{1-4-5} and \eqref{2-3-5} become modular forms which is not easy to observe from their sum-side representations. A central question in this area is to characterize which basic hypergeometric series can be expressed as modular forms.

In this aspect, Nahm \cite{Nahm-1994,Nahm-1995,Nahm-2007} considered a specific class of $q$-hypergeometric series known as Nahm sum or Nahm series
\begin{align*}
f_{A,B,C}(q) := \sum_{\substack{n=(n_1, \dots, n_r)^{\mathrm{T}} \in \mathbb{N}^r}} \frac{q^{\frac{1}{2} n^{\mathrm{T}} A n + n^{\mathrm{T}} B + C}}{(q; q)_{n_1} \cdots (q; q)_{n_r}},
\end{align*}
where $r$ is a positive integer, $A$ is a real positive definite symmetric $r\times r$ matrix, $B$ is an $r$\nobreakdash-di\-men\-sional column vector, and $C$ is a rational scalar. Motivated by considerations from physics, Nahm \cite{Nahm-2007} posed the problem of classifying all rational triples $(A, B, C)$ for which $f_{A,B,C}$ is modular; such a triple $(A,B,C)$ is called as a modular triple.\par
Nahm further conjectured a necessary and sufficient condition on $A$ for it to be the matrix part of a modular triple. In a systematic study, Zagier \cite{Zagier-2007} showed that there are exactly seven modular triples in the rank $r=1$ case:
\begin{align*}
& (1/2,0,-1/40), \quad (1/2,1/2,1/40), \quad (1,0,-1/48), \quad (1,1/2,1/24),\\
& (1,-1/2,1/24),  \quad (2,0,-1/60), \quad (2,1,11/60).
\end{align*}
Notably, the last two triples correspond to the Rogers--Ramanujan identities \eqref{1-4-5} and \eqref{2-3-5}.\par
The following identity is first given by Ramanujan, which is commonly known as Ramanujan's $_1\psi_1$ summation \cite[Appendix~II.29]{red}
\begin{align}\label{1psi1}
 \sum_{k=-\infty}^{\infty} \frac{(a; q)_k}{(b; q)_k} z^k = \frac{(q, az, q/az, b/a; q)_{\infty}}{(b, z, b/az, q/a; q)_{\infty}}, \qquad |b/a| < |z| < 1.
\end{align}
The identity is a bilateral extension of the $q$-binomial theorem \cite[Appendix~II.3]{red}
\begin{align*}
 \sum_{k=0}^{\infty} \frac{(a; q)_k}{(q; q)_k} z^k = \frac{(az; q)_{\infty}}{(z; q)_{\infty}}, \qquad |z| < 1,
\end{align*}
which is a fundamental identity in the theory of basic hypergeometric series. A natural question is whether similar bilateral extensions exist for Rogers--Ramanujan type identities. In 2023, Schlosser \cite{Schlosser-2023} employed bilateral summation formulas to study Nahm sum and derived new bilateral Rogers--Ramanujan type identities by analytic methods.\par
Following the notation of \cite{Wang-2023}, we define a Rogers--Ramanujan type identity of the index $(n_1,n_2,\dots,n_k)$ as a finite sum of the form:
\begin{align*}
\sum_{(i_1,\dots,i_k)\in S}\frac{(-1)^{t(i_1,\dots,i_k)}q^{Q(i_1,\dots,i_k)}}{(q^{n_1};q^{n_1})_{i_1}\cdots(q^{n_k};q^{n_k})_{i_k}}=\prod_{(a,n)\in P}(q^a;q^n)^{r(a,n)}_{\infty},
\end{align*}
 where $t(i_1, {\dots}, i_k)$ is an integer-valued function, $Q(i_1, {\dots}, i_k)$ is a rational polynomial in $i_1,{\dots},i_k$, the $n_1,n_2,\dots,n_k$ are positive integers with $\gcd(n_1,n_2,\dots,n_k)=1$, $S\subset\mathbb{Z}^k$, $P\subset\mathbb{Q}^2$ is finite, and $r(a,n)$ is integer-valued. For instance, the Andrews--Gordon identity \eqref{a-g-identity} and Bressoud's identity \eqref{b-identity} are of indices $(1,1,\dots,1)$ and $(1,1,\dots,1,2)$, respectively.\par
Recently, Cao and Wang \cite{Cao-2023} used the integral method to derive several new Rogers--Ra\-manu\-jan type identities of indices
\begin{align*}
(1,1), \quad (1,2), \quad (1,1,1), \quad (1,1,2), \quad (1,1,3), \quad (1,2,2), \quad (1,2,3), \quad (1,2,4).
\end{align*}
Some of these include additional parameters, thereby giving infinite families of such identities. For example, they proved \cite[Theorem~3.4]{Cao-2023}
\begin{align}\label{cao-34}
\sum_{i,j\geq0}\frac{u^{i-j}q^{\binom{i}{2}+\binom{j+1}{2}+a\binom{j-i}{2}}}{(q;q)_i(q;q)_j} =\frac{\bigl(-uq^a,-q/u,q^{a+1};q^{a+1}\bigr)_{\infty}}{(q;q)_{\infty}}.
\end{align}

Motivated by the integral method and recent work on bilateral Rogers--Ramanujan identities, we present the following main results.
\begin{Theorem}\label{main-theorem}
 For $a\in\mathbb{N^+}$, we have
 \begin{align}\label{main-result-a}
 \sum_{i,j\in \mathbb{Z}}\frac{x^iy^jq^{a\binom{i}{2}+a\binom{j+1}{2}+\binom{j-i}{2}}}{(xq^a;q^a)_i(yq^a;q^a)_j} =\frac{(q;q)_{\infty}}{(xq^a,yq^a;q^a)_{\infty}(q^a;q^a)_{\infty}^{a-1}}\Theta_a\bigl(qy,q^2y,\dots,q^{a-1}y,xy;q^a\bigr),\!\!\!
 \end{align}
 where
\begin{align*}
\Theta_a(x_1,x_2,\dots,x_a;q):=\sum_{n_1,n_2,\dots,n_a\in \mathbb{Z}}q^{n_1^2+\cdots+n_a^2+n_1n_2+n_1n_3+\cdots+n_{a-1}n_a}x_1^{n_1}x_2^{n_2}\cdots x_a^{n_a}.
\end{align*}
In particular, when $a=1$, we have
\begin{align}\label{main-result-1}
 \sum_{i,j \in \mathbb{Z}} \frac{x^i y^j q^{i^2 - ij + j^2}}{(xq; q)_i (yq; q)_j} = \frac{(q; q)_{\infty} \bigl(-xyq, -q/xy, q^2; q^2\bigr)_{\infty}}{(xq, yq; q)_{\infty}}.
\end{align}
When $a=2$, we have
\begin{align}
 &\sum_{i,j\in\mathbb{Z}}\frac{x^iy^jq^{\frac{3}{2}i^2+\frac{3}{2}j^2-ij-\frac{1}{2}i +\frac{1}{2}j}}{\bigl(xq^2;q^2\bigr)_i\bigl(yq^2;q^2\bigr)_j}\nonumber\\
 & \qquad =\frac{(q;q)_{\infty}}{\bigl(xq^2,yq^2,q^2;q^2\bigr)_{\infty}} \Bigl(\Bigl(q^4,-yq^3,\frac{-q}{y};q^4\Bigr)_{\infty}\Bigl(q^{12},-x^2yq^5,\frac{-q^7}{x^2y};q^{12}\Bigr)_{\infty}\nonumber\\
 &\qquad \quad {}+xq\Bigl(q^4,-yq,\frac{-q^3}{y};q^4\Bigr)_{\infty}\Bigl(q^{12},-x^2yq^{11},\frac{-q}{x^2y};q^{12}\Bigr)_{\infty}\Bigr).\label{main-result-2}
\end{align}
Moreover, when $x=y=1$, we have
\begin{align}\label{main-result-3}
\Theta_a\bigl(q,q^2,\dots,q^{a-1},1;q^a\bigr) =\frac{(q^a;q^a)^a_{\infty}\bigl(-q,-q^a,q^{a+1};q^{a+1}\bigr)_{\infty}}{(q;q)_{\infty}}.
\end{align}
\end{Theorem}
As applications, these identities yield the following.
\begin{Corollary}\label{cor}
We have
\begin{align}
& \sum_{i,j\geq0}\frac{q^{i^2-ij+j^2}}{(q;q)_i(q;q)_j} =\frac{\bigl(-q, -q, q^2; q^2\bigr)_{\infty}}{(q; q)_{\infty}},\label{cor11}\\
& \sum_{i,j,k \geq0} \frac{ q^{i^2+j^2+k^2+ik+jk}}{(q; q)_i (q; q)_j(q;q)_k} = \frac{\bigl(-q, -q, q^2; q^2\bigr)_{\infty}}{(q; q)_{\infty}},\label{cor111}\\
& \sum_{n_1,n_2,\dots,n_\ell\geq0}
 \frac{q^{Q(n_1,n_2,\dots,n_\ell)}}
 {(q;q)_{n_1}(q;q)_{n_2}\cdots(q;q)_{n_\ell}} =\frac{\bigl(-q, -q, q^2; q^2\bigr)_{\infty}}{(q; q)_{\infty}}, \label{cor1k}\\
& \sum_{i,j\geq0}\frac{q^{\frac{3}{2}i^2+\frac{3}{2}j^2-ij-\frac{1}{2}i+\frac{1}{2}j}}{\bigl(q^2;q^2\bigr)_i\bigl(q^2;q^2\bigr)_j} =\frac{\bigl(q^3,q^3,q^6;q^6\bigr)_{\infty}}{(q;q)_{\infty}},\label{cor22}\\
& \sum_{i,j\geq0}\frac{q^{2i^2+2j^2-ij-i+j}}{\bigl(q^3,q^3\bigr)_i\bigl(q^3,q^3\bigr)_j} =\frac{\bigl(q^2,q^2\bigr)^2_{\infty}}{(q,q)_{\infty}\bigl(q^3,q^3\bigr)_{\infty}},\label{cor33}\\
& \sum_{i,j\geq0}\frac{q^{3i^2+3j^2-ij-2i+2j}}{\bigl(q^5,q^5\bigr)_i\bigl(q^5,q^5\bigr)_j} =\frac{\bigl(q^2,q^2\bigr)^2_{\infty}\bigl(q^3,q^3\bigr)_{\infty} \bigl(q^{12},q^{12}\bigr)_{\infty}}{(q,q)_{\infty}\bigl(q^4,q^4\bigr)_{\infty} \bigl(q^5,q^5\bigr)_{\infty}\bigl(q^6,q^6\bigr)_{\infty}},\label{cor55}
\end{align}
where $Q(n_1,n_2,\dots,n_\ell)=(n_1+n_3+\cdots+n_\ell)^2-(n_1+n_3+\cdots+n_\ell)(n_2+n_3+\cdots+n_\ell)+(n_2+n_3+\cdots+n_\ell)^2+
 \sum_{i \in \{4,\dots,\ell\}} (n_1+n_i+\cdots+n_\ell)(n_2+n_i+\cdots+n_\ell)+n_1n_2$ with integers $\ell\geq4$.
\end{Corollary}
\begin{Remark}
Setting $u=a=1$ in \eqref{cao-34}, we can also obtain \eqref{cor11}. This identity is a special case of a more general family of identities for Nahm sums associated with the tensor product of Cartan matrices \cite[equation~(1.6)]{Warnaar}. Setting $\nu=0$ in \cite[equation~(4.37)]{Wang-2024}, we can also obtain~\eqref{cor111}. Setting $\alpha=\frac{3}{2}$, $\nu=-\frac{1}{6}$ and substituting $q\rightarrow q^2$ in \cite[equation~(3.1)]{Wang-2024-2}, we can also obtain~\eqref{cor22}. Setting $\alpha=\frac{4}{3}$, $\nu=-\frac{1}{4}$ and substituting $q\rightarrow q^3$ in \cite[equation~(3.1)]{Wang-2024-2}, we~can also obtain \eqref{cor33}. Setting $\alpha=\frac{6}{5}$, $\nu=-\frac{1}{3}$ and substituting $q\rightarrow q^5$ in \cite[equation~(3.1)]{Wang-2024-2}, we~can also obtain \eqref{cor55}.
\end{Remark}

The rest of this paper is organized as follows. In Section~\ref{sec-pre}, we review some fundamental identities, the integral method and modular functions. In Section~\ref{proof}, we prove Theorem \ref{main-theorem} using the integral method and derive several multi-sum Rogers--Ramanujan identities as consequences.

\section{Preliminaries}\label{sec-pre}
\subsection{Some fundamental identities and the integral method}
In this section, we recall several fundamental identities and techniques from the theory of $q$-series that will be used throughout this paper.\par
We have that
\begin{align*}
 \lim_{a\rightarrow \infty}(a;q)_na^{-n}=(-1)^{n}q^{\binom{n}{2}}.
\end{align*}
We also require the Jacobi triple product identity \cite[p.~15]{red}
\begin{align*}
 (q,z,q/z;q)_\infty=\sum_{n=-\infty}^{\infty}(-1)^nq^{\binom{n}{2}}z^n.
\end{align*}
Starting from Ramanujan's $_1\psi_1$ summation formula \eqref{1psi1}, we substitute $z\rightarrow z/a$ and take $a\rightarrow \infty$, yielding
\begin{align}\label{bila-euler}
 \sum_{k\in \mathbb{Z}} \frac{(-z)^kq^{\binom{k}{2}}}{(b;q)_k}=\frac{(q,z,q/z;q)_\infty}{(b,b/z;q)_\infty}.
\end{align}

Two special cases of \eqref{bila-euler} will be particularly useful. First, setting $b=xq$ and replacing $z$ with $xz$, we obtain
\begin{align}\label{circle-1}
 \sum_{i\in \mathbb{Z}} \frac{(-xz)^iq^{\binom{i}{2}}}{(xq;q)_i}=\frac{(q,xz,q/xz;q)_\infty}{(xq,q/z;q)_\infty}.
\end{align}
Second, setting $b=yq$ and replacing $z$ with $yq/z$, we obtain
\begin{align}\label{circle-2}
 \sum_{j\in \mathbb{Z}} \frac{(-yq/z)^jq^{\binom{j}{2}}}{(yq;q)_j}=\frac{(q,yq/z,z/y;q)_\infty}{(yq,z;q)_\infty}.
\end{align}

For a Laurent series \( f(z) = \sum_{n = -\infty}^\infty a(n) z^n \), we denote by \([z^n] f(z)\) the coefficient of \( z^n \), i.e., \([z^n] f(z) = a(n)\). A key observation is that this coefficient can be extracted via contour integration
\begin{align*}
 \oint_K f(z) \frac{{\rm d}z}{2\pi {\rm i} z} = [z^0] f(z),
\end{align*}
where \( K \) is any positively oriented simple closed contour encircling the origin. This principle underlies the integral method, which we employ to prove Rogers--Ramanujan type identities. The general approach is to express the sum side of an identity as a finite linear combination of integrals involving infinite products, and then evaluate each integral explicitly.

\subsection{Modular functions}
According to \cite{Berndt-III} and \cite{Ra77}, let $\mathcal{H}=\{\tau \mid \operatorname{Im}(\tau) > 0\}$ denote the complex upper half-plane. For each \[M=\begin{pmatrix} a & b \\c & d \end{pmatrix} \in M_2^{+}(\mathbb{Z}),\] where $M_2^{+}(\mathbb{Z})$ is the set of integer $2\times 2$ matrices with positive determinant, the bilinear transformation $M(\tau)$ is defined by
\begin{align*}
M\tau=M(\tau)=\frac{a\tau+b}{c\tau+d}.
\end{align*}
The slash operator is defined by
\begin{align*}
(f|M)(\tau)=f(M\tau),
\end{align*}
and satisfies
\begin{align*}
f|MS=f|M|S,
\end{align*}
for matrices $M$ and $S$. The modular group $\Gamma(1)$ is defined by
\begin{align*}
\Gamma(1)=\left\{\begin{pmatrix} a & b \\c & d \end{pmatrix} \in M_2^{+}(\mathbb{Z}) \mid ad-bc=1 \right\}.
\end{align*}
We consider the following subgroup $\Gamma$ of the modular group with finite index
\begin{align*}
\Gamma_1(N)=\left\{\begin{pmatrix} a & b \\c & d \end{pmatrix} \in \Gamma(1) \mid \begin{pmatrix} a & b \\c & d \end{pmatrix}  \equiv \begin{pmatrix} 1 & * \\0 & 1 \end{pmatrix} \ ({\rm mod}\, N) \right\}.
\end{align*}
The Dedekind eta function is defined to be
\begin{align*}
\eta_m(\tau)=\eta(m\tau)=q^{\frac{m}{24}}\prod_{n=1}^{\infty}(1-q^{mn}).
\end{align*}
The generalized Dedekind eta function is defined to be
\begin{align}\label{generalized-eta}
\eta_{\delta,g}(\tau)=q^{\frac{\delta}{2}P_2(g/\delta)}\prod_{m\equiv\pm g\,({\rm mod}\,\delta)}(1-q^m),
\end{align}
where $P_2(t)=\{t\}^2-\{t\}+\frac{1}{6}$ is the second periodic Bernoulli polynomial, $\{t\}=t-[t]$ is the fractional part of $t, g, \delta, m\in\mathbb{Z}^{+}$ and $0 < g < \delta$. The function $\eta_{\delta,g}(\tau)$ is a modular function on ${\rm SL}_2(\mathbb{Z})$ with a multiplier system. Let $N$ be a fixed positive integer. A generalized Dedekind eta-product of level $N$ has the form
\begin{align}\label{f-tau}
f(\tau)=\prod_{\delta|N,~0 < g < \delta}\eta_{\delta,g}^{r_{\delta,g}}(\tau),
\end{align}
where
\begin{align*}
r_{\delta,g}\in \begin{cases}
			\frac{1}{2}\mathbb{Z},  &g=\frac{\delta}{2},\\
			\mathbb{Z}, &\text{otherwise}.\end{cases}
\end{align*}
Robins \cite{Ro94} has found sufficient conditions under which a generalized eta-product is a modular function on $\Gamma_1(N)$.
\begin{Theorem}[{\cite[Theorem 3]{Ro94}}]\label{Thm-2.1}
The function $f(\tau)$, defined in~\eqref{f-tau}, is a modular function on~$\Gamma_1(N)$ if
\begin{alignat*}{3}
&(i) \quad && \sum_{\delta|N,~g}\delta P_2\Bigl(\frac{g}{\delta}\Bigr)r_{\delta,g} \equiv 0 \pmod 2, \quad \text{and}& \\
&(ii) \quad && \sum_{\delta|N,~g}\frac{N}{\delta}P_2(0)r_{\delta,g} \equiv 0 \pmod 2.&
\end{alignat*}
\end{Theorem}
Cho, Koo and Park \cite{CKP09} have found a set of inequivalent cusps for $\Gamma_1(N)\cap\Gamma_0(mN)$. The group $\Gamma_1(N)$ corresponds to the case $m=1$.

\begin{Theorem}[{\cite[Corollary 4]{CKP09}}]\label{Thm-2.2}
Let $a, c, a', c'\in\mathbb{Z}$ with $(a, c)=(a', c')=1$.
\begin{itemize}\itemsep=0pt
 \item[$(i)$] The cusps $\frac{a}{c}$ and $\frac{a'}{c'}$ are equivalent mod $\Gamma_1(N)$ if and only if \[ \binom{a'}{c'}\equiv \pm\binom{a+nc}{c} \pmod N \] for some integers~$n$.
 \item[$(ii)$] The following is a complete set of inequivalent cusps mod $\Gamma_{1}$:
\begin{align*}
\mathcal{S}={} &\biggl\{\frac{y_{c,j}}{x_{c,j}}\,\Bigl|\, 0 < c|N,\, 0 < s_{c,i},a_{c,j}\leq N,\, (s_{c,i},N)=(a_{c,j},N)=1,\\
&s_{c,i}=s_{c,i'}\Leftrightarrow s_{c,1}\equiv\pm s_{c',i'} \ \Big({\rm mod} \ \frac{N}{c}\Big),\\
&a_{c,j}=a_{c,j'}\Leftrightarrow \begin{cases}
a_{c,j}\equiv\pm a_{c,j'} \ ({\rm mod}~c), & \text{if}~c=\frac{N}{2}~\text{or}~ N,\\
a_{c,j}\equiv\pm a_{c,j'}\ ({\rm mod}~c), & \text{otherwise}.\end{cases}\\
&x_{c,i},y_{c,j}\in\mathbb{Z}~\text{chosen so that}~x_{c,i}\equiv cs_{c,i} ,\,  y_{c,j}\equiv a_{c,j}\ ({\rm mod}~N), \, (x_{c,i},y_{c,j})=1\biggr\}.
\end{align*}
 \item[(iii)] The fan width of the cusp $\frac{a}{c}$ is given by
\begin{align*}
\kappa\Bigl(\frac{a}{c},\Gamma_1(N)\Bigr)=\begin{cases}
1, &\text{if}~N=4~\text{and}~(c,4)=2,\\
\dfrac{N}{(c,N)}, &\text{otherwise}.\end{cases}
\end{align*}
\end{itemize}
\end{Theorem}
In this theorem, it is understood, as usual that the fraction $\frac{\pm 1}{0}$ corresponds to $\infty$. Robins~\cite{Ro94} has calculated the invariant order of $\eta_{\delta,g}(\tau)$ at any cusp. This gives a method for calculating the invariant order at any cusp of a generalized eta-product.
\begin{Theorem}[\cite{Ro94}]\label{Thm-2.3}
The order at the cusp $\zeta=\frac{a}{c}$ $($assuming $(a,c)=1)$ of the generalized eta-function $\eta_{\delta,g}(\tau)$ $($defined in \eqref{generalized-eta} and assuming $0 < g < \delta)$ is
\begin{align*}
ord(\eta_{\delta,g}(\tau);\zeta)=\frac{\varepsilon^2}{2\delta}P_2\Bigl(\frac{ag}{\varepsilon}\Bigr),
\end{align*}
where $\varepsilon=(\delta,c)$.
\end{Theorem}
\begin{Theorem}[{\cite[Corollary 2.5]{FrGa19}}]\label{Thm-2.4}
Let $f_1(\tau), f_2(\tau), \dots, f_n(\tau)$ be generalized eta-products that are modular functions on $\Gamma_1(N)$. Let $\mathcal{S}_N$ be a set of inequivalent cusps for~$\Gamma_1(N)$. Define the constant
\begin{align}\label{B}
B=\sum_{s\in\mathcal{S}_N,~s\neq\infty}\min(\{\operatorname{Ord}(f_j,s,\Gamma_1(N)) \mid 1 \leq j\leq n\}\cup\{0\}),
\end{align}
and consider
\begin{align*}
g(\tau):=\alpha_1f_1(\tau)+\alpha_2f_2(\tau)+\cdots+\alpha_nf_n(\tau)+1,
\end{align*}
where each $\alpha_j\in\mathbb{C}$. Then
$g(\tau)=0$
if and only if
$\operatorname{Ord}(g(\tau),\infty,\Gamma_1(N)) > -B$.
\end{Theorem}
A more comprehensive description can be found in \cite{FrGa19}. We have utilized a MAPLE package known as \texttt{thetaids}, which facilitates the implementation of the aforementioned algorithm.\footnote{See
\url{http://qseries/org/fgarvan/qmaple/thetaids/}.}

\section{The proof of Theorem \ref{main-theorem} and Corollary \ref{cor}}\label{proof}
In this section, we prove Theorem \ref{main-theorem} using the integral method. We also show how special choices of parameters in Theorem \ref{main-theorem} lead to several Rogers--Ramanujan type identities.
\begin{proof}[The proof of Theorem \ref{main-theorem}]
First of all, we notice a fact about the $q$-shifted factorials:
\begin{align}\label{key-step}
\frac{(z,q/z;q)_{\infty}}{(z,q^a/z;q^a)_{\infty}}=\prod_{l=1}^{a-1}\bigl(zq^l,q^{a-l}/z;q^a\bigr)_{\infty}.
\end{align}
Then we consider the following contour integral:
\begin{align*}
I:={}&\oint \frac{(q^a,xz,q^a/xz;q^a)_{\infty}(q^a,yq^a/z,z/y;q^a)_{\infty}} {(xq^a,q^a/z;q^a)_{\infty}(yq^a,z;q^a)_{\infty}}(q,z,q/z;q)_{\infty}\frac{{\rm d}z}{2\pi {\rm i}z}\\
={}&\frac{(q;q)_{\infty}}{(xq^a,yq^a;q^a)_{\infty}(q^a;q^a)^{a-1}_{\infty}} \\ & \times\oint(q^a,xz,q^a/xz;q^a)_{\infty}(q^a,yq^a/z,z/y;q^a)_{\infty}
\prod_{l=1}^{a-1}\bigl(q^a,zq^l,q^{a-l}/z;q^a\bigr)_{\infty}\frac{{\rm d}z}{2\pi {\rm i}z}\\
={}& \frac{(q;q)_{\infty}}{(xq^a,yq^a;q^a)_{\infty}(q^a;q^a)^{a-1}_{\infty}}\\
& \times\oint\sum_{i\in\mathbb{Z}}(-1)^iq^{a\binom{i}{2}}(xz)^i\sum_{j\in\mathbb{Z}} (-1)^jq^{a\binom{j}{2}}\Bigl(\frac{yq^a}{z}\Bigr)^j
\prod_{l=1}^{a-1}\sum_{n_l\in\mathbb{Z}}(-1)^{n_l}q^{a\binom{n_l}{2}}\bigl(zq^l\bigr)^{n_l}\frac{{\rm d}z}{2\pi {\rm i}z}.
\end{align*}
Here we use \eqref{key-step} for the second equality and the third equality follows the Jacobi triple product identity.

Multiplying the series and extracting the coefficient of $z^0$, we require $i-j+n_1+\cdots+n_{a-1}=0$. We let $i=n_a$ and then $j=\sum_{l=1}^an_l$. Thus,
\begin{align*}
I&=\frac{(q;q)_{\infty}}{(xq^a,yq^a;q^a)_{\infty}(q^a;q^a)^{a-1}_{\infty}} \!\sum_{n_1,n_2,\dots,n_a\in\mathbb{Z}}\!\!\!\! x^{n_a}y^{\sum_{l=1}^an_l}q^{a\binom{\sum_{l=1}^an_l}{2}+a\sum_{l=1}^an_l+a\sum_{l=1}^{a}\binom{n_l}{2}+\sum_{l=1}^{a-1}ln_l}\\
&=\frac{(q;q)_{\infty}}{(xq^a,yq^a;q^a)_{\infty}(q^a;q^a)_{\infty}^{a-1}}\Theta_a\bigl(qy,q^2y,\dots,q^{a-1}y,xy;q^a\bigr).
\end{align*}
On the other hand, starting from \eqref{circle-1} and \eqref{circle-2}, we substitute $q\rightarrow q^a$. Then the same integral can be written as:
\begin{align*}
I&=\oint\sum_{i\in\mathbb{Z}}\frac{(-xz)^iq^{a\binom{i}{2}}}{(xq^a;q^a)_i}\sum_{j\in\mathbb{Z}} \frac{(-yq^a/z)^jq^{a\binom{j}{2}}}{(xq^a;q^a)_j}\sum_{k\in\mathbb{Z}}(-1)^kq^{\binom{k}{2}}z^k\frac{{\rm d}z}{2\pi {\rm i}z}\\
&=\sum_{i,j\in \mathbb{Z}}\frac{x^iy^jq^{a\binom{i}{2}+a\binom{j+1}{2}+\binom{j-i}{2}}}{(xq^a;q^a)_i(yq^a;q^a)_j}.
\end{align*}
Similarly, the second equality holds because we require $i-j+k=0$, i.e., $k=j-i$ to extract the coefficient of $z^0$. Equating the two expressions for $I$ completes the proof of~\eqref{main-result-a}.

Setting $a=1$ in \eqref{main-result-a}, we can obtain~\eqref{main-result-1} by the Jacobi triple product identity:
\begin{align*}
 \sum_{i,j \in \mathbb{Z}} \frac{x^i y^j q^{i^2 - ij + j^2}}{(xq; q)_i (yq; q)_j}
 & =\frac{(q;q)_\infty}{(xq,yq;q)_\infty}\sum_{i\in \mathbb{Z}}(xy)^iq^{i^2} \\
 & =\frac{(q;q)_\infty}{(xq,yq;q)_\infty}\bigl(-xyq,-q/xy,q^2;q^2\bigr)_\infty.
\end{align*}

Setting $a=2$ in \eqref{main-result-a}, we can obtain
\begin{align}\label{key-2}
\sum_{i,j\in\mathbb{Z}}\frac{x^iy^jq^{\frac{3}{2}i^2+\frac{3}{2}j^2-ij-\frac{1}{2}i +\frac{1}{2}j}}{\bigl(xq^2;q^2\bigr)_i\bigl(yq^2;q^2\bigr)_j} =\frac{(q;q)_{\infty}}{\bigl(xq^2,yq^2,q^2;q^2\bigr)_{\infty}}\sum_{i,j\in\mathbb{Z}}x^jy^{i+j}q^{2i^2+2j^2+2ij+i}.
\end{align}
Then by the Jacobi triple product identity, we have
\begin{align*}
&\sum_{i,j\in\mathbb{Z}}x^jy^{i+j}q^{2i^2+2j^2+2ij+i}\\
& \quad =\sum_{i,j\in\mathbb{Z}}x^jy^{i+j}q^{2[(i+\frac{1}{2}j)^2+\frac{3}{4}j^2]+i}\\
& \quad = \sum_{i,j\in\mathbb{Z}}x^{2j}y^{i+2j}q^{2[(i+j)^2+3j^2]+i}+\sum_{i,j\in\mathbb{Z}}x^{2j-1}y^{i+2j-1}q^{2[(i+j-\frac{1}{2})^2+3j^2-3j+\frac{3}{4}]+i}\\
& \quad =\sum_{k,j\in\mathbb{Z}}x^{2j}y^{k+j}q^{2(k^2+3j^2)+k-j}+\sum_{k,j\in\mathbb{Z}}x^{2j-1}y^{k+j-1}q^{2(k^2-k+3j^2-3j+1)+k-j}\\
& \quad =\sum_{k\in\mathbb{Z}}y^kq^{2k^2+k}\sum_{j\in\mathbb{Z}}(x^2y)^jq^{6j^2-j}+xq\sum_{k\in\mathbb{Z}}y^kq^{2k^2-k}\sum_{j\in\mathbb{Z}}(x^2y)^jq^{6j^2+5j}\\
& \quad =\Bigl(q^4,-yq^3,\frac{-q}{y};q^4\Bigr)_{\infty}\Bigl(q^{12},-x^2yq^5,\frac{-q^7}{x^2y};q^{12}\Bigr)_{\infty}\\
&\qquad +xq\Bigl(q^4,-yq,\frac{-q^3}{y};q^4\Bigr)_{\infty}\Bigl(q^{12},-x^2yq^{11},\frac{-q}{x^2y};q^{12}\Bigr)_{\infty}.
\end{align*}
Combining the above identity and \eqref{key-2}, we complete the proof of \eqref{main-result-2}.

From Wang \cite[equation~(3.1)]{Wang-2024-2}, we have
\begin{align}\label{wang-3.1}
\sum_{i,j\geq0}\frac{q^{\frac{\alpha}{2}i^2+(1-\alpha)ij+\frac{\alpha}{2}j^2+\alpha\nu i-\alpha\nu j}}{(q;q)_i(q;q)_j} =\frac{\bigl(-q^{\frac{\alpha}{2}+\alpha\nu},-q^{\frac{\alpha}{2}-\alpha\nu},q^{\alpha};q^{\alpha}\bigr)_{\infty}}{(q;q)_{\infty}}.
\end{align}
Setting $\alpha=\frac{a+1}{a},~\nu=-\frac{a-1}{2(a+1)}$ and substituting $q\rightarrow q^a$ in \eqref{wang-3.1}, we obtain
\begin{align*}
\sum_{i,j\geq0}\frac{q^{\frac{a+1}{2}i^2-ij+\frac{a+1}{2}j^2-\frac{a-1}{2}i+\frac{a-1}{2} j}}{(q^a;q^a)_i(q^a;q^a)_j}=\frac{\bigl(-q,-q^a,q^{a+1};q^{a+1}\bigr)_{\infty}}{(q^a;q^a)_{\infty}}.
\end{align*}
Setting $x=y=1$ in \eqref{main-result-a}, we obtain
\begin{align*}
\sum_{i,j\geq0}\frac{q^{\frac{a+1}{2}i^2-ij+\frac{a+1}{2}j^2-\frac{a-1}{2}i+\frac{a-1}{2} j}}{(q^a;q^a)_i(q^a;q^a)_j}=\frac{(q;q)_{\infty}}{(q^a;q^a)^{a+1}_{\infty}}\Theta_a\bigl(q,q^2,\dots,q^{a-1},1;q^a\bigr).
\end{align*}
Combining the above two identities, we can prove \eqref{main-result-3}.
\end{proof}

\begin{proof}[The proof of Corollary \ref{cor}]
Setting $x=y=1$ in \eqref{main-result-1} and noting that $1/(q;q)_{-n}=0$ for $n>0$, the bilateral sum reduces to a sum over non-negative indices, yielding \eqref{cor11}.

To prove \eqref{cor1k}, we use the identity \cite[p.~20]{Andrews-1981}
\begin{align*}
\frac{1}{(q;q)_i(q;q)_j}=\sum_{k\geq0}\frac{q^{(i-k)(j-k)}}{(q;q)_k(q;q)_{i-k}(q;q)_{j-k}}.
\end{align*}
Substituting this into~\eqref{cor11} and shifting the summation indices by setting $i\rightarrow i+k$ and $j\rightarrow j+k$, we obtain identity~\eqref{cor111}. Finally, applying the same identity from \cite[p.~20]{Andrews-1981} iteratively yields identity~\eqref{cor1k}.

Setting $x=y=1$ in \eqref{main-result-2}, we obtain
\begin{align*}
\sum_{i,j\geq0}\frac{q^{\frac{3}{2}i^2+\frac{3}{2}j^2-ij-\frac{1}{2}i+\frac{1}{2}j}} {\bigl(q^2;q^2\bigr)_i\bigl(q^2;q^2\bigr)_j}=\frac{1}{\bigl(q^2;q^2\bigr)_{\infty}}
\bigl(\bigl(-q^5,-q^7,q^{12};q^{12}\bigr)_{\infty}+q\bigl(-q,-q^{11},q^{12};q^{12}\bigr)_{\infty}\bigr).
\end{align*}
Therefore, \eqref{cor22} holds if we can prove
\begin{align}\label{key-id}
\bigl(-q^5,-q^7,q^{12};q^{12}\bigr)_{\infty}+q\bigl(-q,-q^{11},q^{12};q^{12}\bigr)_{\infty}
=\frac{\bigl(q^3;q^3\bigr)^2_{\infty}\bigl(q^2;q^2\bigr)_{\infty}}{\bigl(q^6;q^6\bigr)_{\infty}(q;q)_{\infty}}.
\end{align}

We first rewrite \eqref{key-id} as the following modular function identity for generalized eta-products on $\Gamma_1(24)$:
\begin{align}\label{getaid}
0=1-\frac{\eta_{24,1}\eta_{24,10}\eta_{24,11}}{\eta_{24,3}\eta_{24,6}\eta_{24,9}}-\frac{\eta_{24,2}\eta_{24,5}\eta_{24,7}}{\eta_{24,3}\eta_{24,6}\eta_{24,9}}.
\end{align}
We use Theorem~\ref{Thm-2.1} to check that each generalized eta-product is a modular function on~$\Gamma_1(24)$. We then use Theorems~\ref{Thm-2.2} and~\ref{Thm-2.3} to calculate the order of each generalized eta-product at each cusp of~$\Gamma_1(24)$. We calculate the constant in equation~\eqref{B} to find that $B=-4$. We let $g(\tau)$ be the right-hand side of~\eqref{getaid} and easily show that $\operatorname{Ord}(g(\tau),\infty,\Gamma_1(24)) > 4$. Thus~\eqref{getaid} follows by Theorem~\ref{Thm-2.4}.

Setting $x=y=1$ and $a=3$ in \eqref{main-result-a}, we obtain
\begin{align*}
\sum_{i,j\geq0}\frac{q^{2i^2+2j^2-ij-i+j}}{\bigl(q^3,q^3\bigr)_i\bigl(q^3,q^3\bigr)_j} =\frac{(q,q)_{\infty}}{\bigl(q^3,q^3\bigr)^4_{\infty}}\sum_{i,j,k\in\mathbb{Z}}q^{3(i^2+j^2+k^2+ij+ik+jk)+i+2j}.
\end{align*}
Therefore, \eqref{cor33} holds if we can prove
\begin{align}\label{key-33}
\sum_{i,j,k\in\mathbb{Z}}q^{3(i^2+j^2+k^2+ij+ik+jk)+i+2j}=\frac{\bigl(q^3,q^3\bigr)^3_{\infty} \bigl(q^2,q^2\bigr)^2_{\infty}}{(q,q)^2_{\infty}}.
\end{align}
Firstly, we have
\begin{align*}
&\sum_{i,j,k\in\mathbb{Z}}q^{3(i^2+j^2+k^2+ij+ik+jk)+i+2j}=\sum_{i,j,k\in\mathbb{Z}}q^{3[(i+\frac{j+k}{2})^2+(\frac{j+k}{2})^2+\frac{j^2+k^2}{2}]+i+2j}\\
& \qquad =\sum_{i,j,k\in\mathbb{Z}}q^{3[(i+j+\frac{k+1}{2})^2+(j+\frac{k+1}{2})^2+2j^2+2j+\frac{k^2+1}{2}]+i+4j+2} \\ & \qquad +\sum_{i,j,k\in\mathbb{Z}}q^{3[(i+j+\frac{k}{2})^2+(j+\frac{k}{2})^2+2j^2+\frac{k^2}{2}]+i+4j}\\
& \qquad =\sum_{i,j,k\in\mathbb{Z}}q^{3[(i+j+k)^2+(j+k)^2+2j^2+2k^2]+i+4j} \\
& \qquad\quad{} +\sum_{i,j,k\in\mathbb{Z}}q^{3[(i+j+k+\frac{1}{2})^2+(j+k+\frac{1}{2})^2+2j^2+2k^2+2k+\frac{1}{2}]+i+4j}\\
&\qquad\quad{} +\sum_{i,j,k\in\mathbb{Z}}q^{3[(i+j+k+\frac{1}{2})^2+(j+k+\frac{1}{2})^2+2j^2+2j+2k^2\frac{1}{2}]+i+4j+2}\\
&\qquad\quad{} +\sum_{i,j,k\in\mathbb{Z}}q^{3[(i+j+k)^2+(j+k)^2+2j^2+2j+2k^2-2k+1]+i+4j+2}.
\end{align*}
The method to simplify these four series is similar, so we take the first series as an example. Setting $j+k=m$ and $i+m=l$, we obtain
\begin{align*}
&\sum_{i,j,k\in\mathbb{Z}}q^{3[(i+j+k)^2+(j+k)^2+2j^2+2k^2]+i+4j}\\
& \qquad =\sum_{l,m,k\in\mathbb{Z}}q^{3(m^2+l^2+2(m-k)^2+2k^2)+l+3m-4k} \\
& \qquad =\sum_{l,m,k\in\mathbb{Z}}q^{3(2m^2+l^2+(m-2k)^2)+l+3m-4k}\\
& \qquad =\sum_{l,m,k\in\mathbb{Z}}q^{3(8m^2+l^2+(2m-2k)^2)+l+6m-4k}+\sum_{l,m,k\in\mathbb{Z}}q^{3(2(2m+1)^2+l^2+(2m-2k+1)^2)+l+6m+3-4k}\\
& \qquad =\sum_{l,m,n\in\mathbb{Z}}q^{3(8m^2+l^2+4n^2)+l+2m+4n}+\sum_{l,m,k\in\mathbb{Z}}q^{3(2(2m+1)^2+l^2+(2n+1)^2)+l+2m+3+4n}\\
& \qquad =\sum_{l\in\mathbb{Z}}q^{3l^2+l}\sum_{m\in\mathbb{Z}}q^{24m^2+2m}\sum_{n\in\mathbb{Z}}q^{12n^2+4n}+q^{12}\sum_{l\in\mathbb{Z}}q^{3l^2+l}\sum_{m\in\mathbb{Z}}q^{24m^2+26m}\sum_{n\in\mathbb{Z}}q^{12n^2+16n}\\
& \qquad =\bigl(q^6,-q^2,-q^4;q^6\bigr)_{\infty} \bigl(q^{48},-q^{26},-q^{22};q^{48}\bigr)_{\infty}\bigl(q^{24},-q^{16},-q^8;q^{24}\bigr)_{\infty}\\
&\qquad\quad{} +q^6\bigl(q^6,-q^2,-q^4;q^6\bigr)_{\infty} \bigl(q^{48},-q^2,-q^{46};q^{48}\bigr)_{\infty}\bigl(q^{24},-q^4,-q^{20};q^{24}\bigr)_{\infty}.
\end{align*}
Here we let $m-k=n$ for the fourth equality and the last equality follows the Jacobi triple product identity. Therefore, \eqref{key-33} holds if we can prove
\begin{align}\notag
&\frac{\eta_4\eta^2_6\eta_{16}\eta^2_{24}\eta_{96,44}}{\eta_2\eta_8\eta_{12}\eta_{96,22}\eta_{96,26}}+\frac{\eta^2_6\eta^2_8\eta^2_{48}\eta_{96,4}}{\eta_2\eta_{16}\eta_{24}\eta_{96,2}\eta_{96,46}}+\frac{\eta^2_2
\eta_3\eta_{12}\eta_{24}\eta_{48}\eta_{48,20}\eta_{96,20}}{\eta_1\eta_4\eta_6\eta_{48,10}\eta_{48,14}\eta_{96,10}\eta_{96,38}}\\
&\qquad{} +\frac{\eta^2_2\eta_3\eta_{12}\eta_{24}\eta_{48}\eta_{48,4}\eta_{96,28}}{\eta_1\eta_4\eta_6\eta_{48,2}\eta_{48,22}\eta_{96,14}\eta_{96,34}}+\frac{\eta^2_6\eta^2_8\eta^2_{48}\eta_{96,44}}{\eta_2\eta_{16}\eta_{24}\eta_{96,22}\eta_{96,26}}+
\frac{\eta_4\eta^2_6\eta_{16}\eta^2_{24}\eta_{96,4}}{\eta_2\eta_8\eta_{12}\eta_{96,2}\eta_{96,46}}\notag\\
&\qquad{} +\frac{\eta^2_2\eta_3\eta_{12}\eta_{24}\eta_{48}\eta_{48,4}\eta_{96,20}}
{\eta_1\eta_4\eta_6\eta_{48,2}\eta_{48,22}\eta_{96,10}\eta_{96,38}} +\frac{\eta^2_2\eta_3\eta_{12}\eta_{24}\eta_{48}\eta_{48,20}\eta_{96,28}} {\eta_1\eta_4\eta_6\eta_{48,10}\eta_{48,14}\eta_{96,14}\eta_{96,34}}=\frac{\eta^3_3\eta^2_2}{\eta^2_1}.\label{key-id2}
\end{align}

We first rewrite \eqref{key-id2} as the following modular function identity for generalized eta-products on $\Gamma_1(96)$:
\begin{align}\label{getaid2}
0=1-x(L_1+L_2+L_3+L_4+L_5+L_6+L_7+L_8),
\end{align}
where
\begin{gather*}
x=\frac{\eta_{96,1}\eta_{96,5}\eta_{96,7}\eta^2_{96,11}\eta^2_{96,13}\eta^2_{96,17}\eta^2_{96,19}\eta^2_{96,23}\eta^2_{96,25}\eta^2_{96,29}\eta^2_{96,31}\eta^2_{96,35}\eta^2_{96,37}\eta^2_{96,41}\eta^2_{96,43}\eta^2_{96,47}}
{\eta_{96,3}\eta^3_{96,6}\eta_{96,9}\eta^3_{96,12}\eta_{96,15}\eta^3_{96,18}\eta_{96,21}\eta^3_{96,24}\eta_{96,27}\eta^3_{96,30}\eta_{96,33}\eta^3_{96,36}\eta_{96,39}\eta^3_{96,42}\eta_{96,45}},\\
L_1=\frac{\eta_{96,1}\eta_{96,5}\eta_{96,6}\eta_{96,7}\eta_{96,12}\eta_{96,18}\eta^2_{96,24}\eta_{96,30}\eta_{96,36}\eta_{96,42}\eta_{96,44}}{\eta_{96,2}\eta_{96,8}\eta_{96,10}\eta_{96,14}\eta^2_{96,22}\eta^2_{96,26}\eta_{96,34}\eta_{96,38}\eta_{96,40}\eta_{96,46}},\\
L_2=\frac{\eta_{96,1}\eta_{96,5}\eta_{96,6}\eta_{96,7}\eta_{96,8}\eta_{96,12}\eta_{96,18}\eta^2_{96,24}\eta_{96,30}\eta_{96,36}\eta_{96,40}\eta_{96,42}}{\eta^2_{96,2}\eta_{96,10}\eta_{96,14}\eta_{96,20}\eta_{96,22}\eta_{96,26}\eta_{96,28}\eta_{96,34}\eta_{96,38}\eta_{96,44}\eta^2_{96,46}},\\
L_3=\frac{\eta_{96,2}\eta_{96,6}\eta_{96,12}\eta_{96,18}\eta^2_{96,20}\eta_{96,22}\eta^2_{96,24}\eta_{96,26}\eta_{96,28}\eta_{96,30}\eta_{96,36}\eta_{96,42}\eta_{96,46}}
{\eta_{96,10}\eta_{96,11}\eta_{96,13}\eta_{96,17}\eta_{96,19}\eta_{96,23}\eta_{96,25}\eta_{96,29}\eta_{96,31}\eta_{96,35}\eta_{96,37}\eta_{96,38}\eta_{96,41}\eta_{96,43}\eta_{96,47}},\\
L_4=\frac{\eta_{96,4}\eta_{96,6}\eta_{96,10}\eta_{96,12}\eta_{96,18}\eta^2_{96,24}\eta_{96,28}\eta_{96,30}\eta_{96,36}\eta_{96,38}\eta_{96,42}\eta_{96,44}}
{\eta_{96,11}\eta_{96,13}\eta_{96,17}\eta_{96,19}\eta_{96,23}\eta_{96,25}\eta_{96,29}\eta_{96,31}\eta_{96,35}\eta_{96,37}\eta_{96,41}\eta_{96,43}\eta_{96,47}},\\
L_5=\frac{\eta_{96,4}\eta_{96,6}\eta_{96,12}\eta_{96,14}\eta_{96,18}\eta_{96,20}\eta^2_{96,24}\eta_{96,30}\eta_{96,34}\eta_{96,36}\eta_{96,42}\eta_{96,44}}
{\eta_{96,11}\eta_{96,13}\eta_{96,17}\eta_{96,19}\eta_{96,23}\eta_{96,25}\eta_{96,29}\eta_{96,31}\eta_{96,35}\eta_{96,37}\eta_{96,41}\eta_{96,43}\eta_{96,47}},\\
L_6=\frac{\eta_{96,2}\eta_{96,6}\eta_{96,12}\eta_{96,18}\eta_{96,20}\eta_{96,22}\eta^2_{96,24}\eta_{96,26}\eta^2_{96,28}\eta_{96,30}\eta_{96,36}\eta_{96,42}\eta_{96,46}}
{\eta_{96,11}\eta_{96,13}\eta_{96,14}\eta_{96,17}\eta_{96,19}\eta_{96,23}\eta_{96,25}\eta_{96,29}\eta_{96,31}\eta_{96,34}\eta_{96,35}\eta_{96,37}\eta_{96,41}\eta_{96,43}\eta_{96,47}},\\
L_7=\frac{\eta_{96,1}\eta_{96,5}\eta_{96,6}\eta_{96,7}\eta_{96,8}\eta_{96,12}\eta_{96,18}\eta^2_{96,24}\eta_{96,30}\eta_{96,36}\eta_{96,40}\eta_{96,42}}{\eta_{96,2}\eta_{96,4}\eta_{96,10}\eta_{96,14}\eta_{96,20}\eta^2_{96,22}\eta^2_{96,26}\eta_{96,28}\eta_{96,34}\eta_{96,38}\eta_{96,46}},\\
L_8=\frac{\eta_{96,1}\eta_{96,4}\eta_{96,5}\eta_{96,6}\eta_{96,7}\eta_{96,12}\eta_{96,18}\eta^2_{96,24}\eta_{96,30}\eta_{96,36}\eta_{96,42}}{\eta^2_{96,2}\eta_{96,8}\eta_{96,10}\eta_{96,14}\eta_{96,22}\eta_{96,26}\eta_{96,34}\eta_{96,38}\eta_{96,40}\eta^2_{96,46}}.
\end{gather*}
We use Theorem \ref{Thm-2.1} to check that each generalized eta-product is a modular function on $\Gamma_1(96)$. We then use Theorems~\ref{Thm-2.2} and~\ref{Thm-2.3} to calculate the order of each generalized eta-product at each cusp of $\Gamma_1(96)$. We calculate the constant in equation \eqref{B} to find that $B=-192$. We let~$g(\tau)$ be the right-hand side of \eqref{getaid2} and easily show that $\operatorname{Ord}(g(\tau),\infty,\Gamma_1(96)) > 192$. Thus, \eqref{getaid2} follows by Theorem \ref{Thm-2.4}.

Setting $x=y=1$ and $a=5$ in \eqref{main-result-a}, together with \eqref{main-result-3}, we can prove \eqref{cor55}.
\end{proof}

\subsection*{Acknowledgements}
We are grateful to Rong Chen providing the conjecture about \eqref{main-result-a} and \eqref{main-result-1}. We are grateful to Professor Liuquan Wang for drawing our attention to the references \cite{Wang-2024, Wang-2024-2}, and to Professor Ole Warnaar for bringing the paper \cite{Warnaar} to our notice. We are also grateful to the referee for valuable comments and suggestions. The first author was supported by the National Key R\&D Program of China (Grant No. 2024YFA1014500) and the National Natural Science Foundation of China (Grant No. 12201387).

\pdfbookmark[1]{References}{ref}
\LastPageEnding

\end{document}